\documentclass[twoside]{article}

\usepackage{amsfonts,amsmath,amsthm}
\usepackage[latin1]{inputenc}
\newtheorem{proposition}{Proposition}

\newtheorem{rem}{Remark}
\newtheorem{corollary}{Corollary}

\def\qed{\hbox to 0pt{}\hfill$\rlap{$\sqcap$}\sqcup$}

\newtheorem{thm}{Theorem}
\usepackage{epsfig}
\newcounter{subfigure}

\makeatletter

\begin{document}

\vskip-45mm
\title{On the controllability of some steady states in the case of nonlinear discrete dynamical systems with control}
\author{E. Kaslik$^{1,3}$, A.M. Balint$^{2}$, A. Grigis $^{3}$, St. Balint$^{1}$ }

\date{}
\maketitle

\noindent$^1$Department of Mathematics, West University of
Timi\c{s}oara, Romania \hfill\break e-mail: balint@balint.uvt.ro

\noindent$^2$Department of Physics, West University of
Timi\c{s}oara, Romania

\noindent$^3$L.A.G.A, Institut Galilee, Universite Paris 13,
France

\vspace{1mm}

\begin{abstract}
The main objective of this paper is to show that two
asymptotically stable steady states which belong to an analytic
path of asymptotically stable steady states can be gradually
transferred one to the other by successive changes of the control
parameters.
\end{abstract}

%\footnotetext{\textit{AMS Subject Classification:} 34.K.20 \\
%Keywords:difference equations, fixed point, asymptotically stable,
%domain of attraction}

%%%%%%%%%%%%%%%%%%%%%%%%%%%%%%%%%%%%%%%%%%%%%%%%%%%%%%%%%%%%%%%%%%%%%%%%%%%%%%%%%%%%

\section{Introduction}

For nonlinear systems of differential equations with control, it
has been proved (see \cite{Balint}), that two asymptotically
stable steady states belonging to an analytic path of
asymptotically stable steady states can be transferred one in the
other by successive maneuvers along the path. That means,
according to \cite{Rasvan}, that the process described by the
system can be piloted through the domains of attraction of the
intermediary steady states, from the first to the second steady
state.

In this paper, a similar result is established for the nonlinear
discrete dynamical systems with control. A theorem from
\cite{EBBB} is used, which states that the domain of attraction of
an asymptotically stable fixed point of a nonlinear system of
analytic difference equations, is the natural domain of
analyticity of a certain Lyapunov function.

\section{Preliminaries}

We consider the following nonlinear discrete dynamical system with
control:

\begin{equation}
\label{ec} x_{k+1}=f(x_{k},\alpha)\qquad k=0,1,2...
\end{equation}
In (\ref{ec}), $f$ is a given function $f:\Omega\times
D\rightarrow\mathbb{R}^{n}$, $\Omega\subset\mathbb{R}^{n}$,
$D\subset\mathbb{R}^{m}$ are domains, $x\in\Omega$ is the state
parameter, and $\alpha\in D$ is the control parameter. What
concerns the regularity of $f$, we assume that $f$ is an analytic
function.

A state $x^{0}\in\Omega$ is a \emph{steady state} for (\ref{ec})
if there exists $\alpha_{x^{0}}\in D$ such that
\begin{equation}\label{fix}
    x^{0}=f(x^{0},\alpha_{x^{0}})
\end{equation}

The steady state $x^{0}\in\Omega$ of (\ref{ec}) is "\emph{stable}"
provided that given any ball
$B(x^{0},\varepsilon)=\{x\in\Omega/\|x-x^{0}\|<\varepsilon\}$,
there is a ball
$B(x^{0},\delta)=\{x\in\Omega/\|x-x^{0}\|<\delta\}$ such that if
$x\in B(x^{0},\delta)$ then $f^{k}(x,\alpha_{x^{0}})\in
B(x^{0},\varepsilon)$, for $k=0,1,2,...$ \cite{Kel-Pet}.

If in addition there is a ball $B(x^{0},r)$ such that
$f^{k}(x,\alpha_{x^{0}})\rightarrow x^{0}$ as $k\rightarrow\infty$
for all $x\in B(x^{0},r)$ then the steady state $x^{0}$ is
"\emph{asymptotically stable}" \cite{Kel-Pet}.

The \emph{domain of attraction} $DA(x^{0})$ of the asymptotically
stable steady state $x^{0}$ is the set of initial states $x\in
\Omega$ from which the system converges to the steady state itself
i.e.

\begin{equation}\label{da}
    DA(x^{0})=\{x\in\Omega | f^{k}(x,\alpha_{x^{0}})\stackrel{k\rightarrow\infty}{\longrightarrow}x^{0}\}
\end{equation}

An \emph{analytic path of steady states} of (\ref{ec}) is an
analytic function $\varphi :D_{1}\subset D\rightarrow \Omega$
which satisfies

\begin{equation}\label{path}
    \varphi (\alpha)=f(\varphi (\alpha),\alpha), \qquad \textrm{for any
    }\alpha\in D_{1}
\end{equation}

An \emph{analytic path of asymptotically stable steady states} of
(\ref{ec}) is an analytic path of steady states which are all
asymptotically stable.

A change of control parameters from $\alpha'$ to $\alpha''$ in
(\ref{ec}) is called \emph{maneuver} and is denoted
$\alpha'\rightarrow\alpha''$. The maneuver
$\alpha'\rightarrow\alpha''$ is \emph{successful} on the path
$\varphi$ if $\alpha',\alpha''\in D_{1}$ and the sequence defined
by
\begin{equation}\label{successful}
   x_{k+1}=f(x_{k},\alpha''), \textrm{ }x_{0}=\varphi(\alpha')
\end{equation}
tends to $\varphi(\alpha'')$ as $k\rightarrow\infty$.

The following proposition from \cite{EBBB} concerning the discrete
dynamical systems without control parameters is helpful.

\begin{proposition}
If the analytic function $g:\Delta\rightarrow\Delta$ from the
system
\begin{equation}\label{ec0}
    y_{k+1}=g(y_{k}), \qquad k=1,2,...
\end{equation}
satisfies the following conditions:
\begin{equation}
  g(0) = 0
\end{equation}
\begin{equation}
  \|\partial_{0}g\| < 1
\end{equation}
then $0$ is an asymptotically stable steady state of (\ref{ec0}).
$DA(0)$ is an open subset of $\Delta$ and coincides with the
natural domain of analyticity of the unique solution $V$ of the
iterative first order functional equation
\begin{equation}
\label{ecVsimplu}
\begin{array}{ll}
\left\{\begin{array}{l}
V(g(y))-V(y)=-\|y\|^{2}\\
V(0)=0
\end{array}\right.
\end{array}
\end{equation}
The function $V$ is positive on $DA(0)$ and
$V(y)\stackrel{y\rightarrow y^{0}}{\longrightarrow}+\infty$, for
any $y^{0}\in FrDA(0)$ ($FrDA(0)$ denotes the boundary of
$DA(0)$).
\end{proposition}

\begin{rem}
If $\varphi:D_{1}\subset D\rightarrow \Omega$ is an analytical
path of steady states of (\ref{ec}), then the function
$g(\cdot,\alpha):\Omega-\varphi(\alpha)\rightarrow\Omega-\varphi(\alpha)$,
given by
\begin{equation}\label{g}
    g(y,\alpha)=f(y+\varphi(\alpha),\alpha)-\varphi(\alpha) \qquad \textrm{for } y\in\Omega-\varphi(\alpha)
\end{equation}
is analytic and satisfies $g(0,\alpha) = 0$, for any $\alpha\in
D_{1}$. Therefore, $y=0$ is a steady state of the system
\begin{equation}
\label{ec0'} y_{k+1}=g(y_{k},\alpha) \qquad k=0,1,2... \qquad
\textrm{for any } \alpha\in D_{1}.
\end{equation}

The steady state $x=\varphi(\alpha)$ of (\ref{ec}) is
asymptotically stable if and only if the steady state $y=0$ of the
system (\ref{ec0'}) is asymptotically stable. The relationship
between the domain of attraction of the steady state
$x=\varphi(\alpha)$ of (\ref{ec}) and that of the steady state
$y=0$ of (\ref{ec0'}) is
$DA(\varphi(\alpha))=\varphi(\alpha)+DA(0)$.
\end{rem}

\section{Theoretical results}

We now state an existence theorem for an analytic path of
asymptotically stable steady states of (\ref{ec}).

\begin{thm}
If the analytic function $f$ from (\ref{ec}) satisfies:
\begin{enumerate}
    \item there exist $(x^{0},\alpha^{0})\in\Omega\times D$ such
    that $x^{0}=f(x^{0},\alpha^{0})$
    \item $\|\partial_{x^{0}}f(x^{0},\alpha^{0}) \|<1$
\end{enumerate}
then there exists a maximal domain $D_{1}\subset D$ containing
$\alpha^{0}$ and a unique analytic path $\varphi :D_{1}\rightarrow
\Omega$ of asymptotically stable steady states of (\ref{ec})
satisfying the following conditions:
\begin{itemize}
    \item[a.] $\varphi(\alpha^{0})=x^{0}$;
    \item[b.] $\|\partial_{x}f(\varphi(\alpha),\alpha)\|<1$ for
    any $\alpha\in D_{1}$;
    \item[c.] For $\alpha',\alpha''\in D_{1}$ the maneuver $\alpha'\rightarrow \alpha''$ is successful on the
    branch $\varphi$ if and only if $\varphi(\alpha')$ belongs to
    the domain of attraction of $\varphi(\alpha'')$.
\end{itemize}
\end{thm}

\begin{proof}
    As the functions $f$ and $\partial_{x}f$ are continuous on $\Omega\times
    D$, taking into account the properties 1. and 2. of $f$, there
    exist two maximal domains $\Omega_{1}\subset\Omega$ and
    $D_{1}$ and a unique analytic function $\varphi :D_{1}\rightarrow\Omega_{1}$ such that:
    \begin{enumerate}
        \item $(x^{0},\alpha^{0})\in\Omega_{1}\times D_{1}$ and $\varphi(\alpha^{0})=x^{0}$;
        \item $\varphi (\alpha)=f(\varphi (\alpha),\alpha)$, for any $\alpha\in
        D_{1}$;
        \item $\|\partial_{x}f(\varphi(\alpha),\alpha)\|<1$ for any $\alpha\in
        D_{1}$.
    \end{enumerate}
This means that $\varphi$ is path of asymptotically stable steady
states for (\ref{ec}) (see Proposition 1 and Remark 1).

A maneuver $\alpha'\rightarrow\alpha''$ is successful on the path
$\varphi$ if and only if the sequence given by (\ref{successful})
tends to $\varphi(\alpha'')$ as $k\rightarrow\infty$, which means
that $\varphi(\alpha')$ belongs to the domain of attraction of
$\varphi(\alpha'')$.
\end{proof}

In the followings, it is assumed that the conditions of Theorem 1
are fulfilled and thus, there exists an analytic path $\varphi$ of
asymptotically stable steady states of (\ref{ec}).

\begin{thm}
    Let be $\varphi:D_{1}\rightarrow \Omega$ an analytic path of asymptotically
    stable steady states of (\ref{ec}). There exist an open set $G\subset\Omega\times D$ and a
    non-negative analytic function $V$ defined on $G$ satisfying the following conditions:
\begin{itemize}
    \item[a.] $G\supset\Gamma=\{(\varphi(\alpha),\alpha)/\alpha\in D_{1}\}$
    \item[b.] \begin{equation}\label{probLyap'}
                \left\{
                \begin{array}{ll}
                 V(f(x,\alpha),\alpha)-V(x,\alpha)=-\|x-\varphi(\alpha)\|^2 \\
                 V(\varphi(\alpha),\alpha)=0 \\
                \end{array}\right.
              \end{equation}
    \item[c.] For any $\alpha\in D_{1}$, $DA(\varphi(\alpha))$ is
    the natural domain of analyticity of $x\rightarrow V(x,\alpha)$
    \item[d.] $V(x,\alpha)\stackrel{x\rightarrow x^{0}}{\longrightarrow}+\infty$,
     for any $x^{0}\in FrDA(\varphi(\alpha))$.
\end{itemize}
\end{thm}

\begin{proof}
Let be $G=\bigcup\limits_{\alpha\in D_{1}}
(DA(\varphi(\alpha))\times \{\alpha\}) \subset \Omega\times D_{1}$
and $V:G\rightarrow \mathbb{R}^{1}_{+}$ defined by
\begin{equation}\label{fctLyap}
    V(x,\alpha)=\sum_{k=0}^{\infty}\|f^{k}(x,\alpha)-\varphi(\alpha)\|^2
\end{equation}

Proposition 1 and Remark 1 provide that the set $G$ and the
function $V(x,\alpha)$ satisfy the conditions a-d.
\end{proof}

\begin{corollary}
    If $\varphi:D_{1}\rightarrow \Omega$ is an analytic path of asymptotically
    stable steady states of (\ref{ec}) then
    for any $\alpha\in D_{1}$ there is an open neighborhood
    $U_{\alpha}$ of $\alpha$ and an open neighborhood $U_{\varphi(\alpha)}$ of $\varphi(\alpha)$
    such that:
    \begin{itemize}
        \item[1.] $\varphi(\alpha')\in U_{\varphi(\alpha)}$, for
        any $\alpha'\in U_{\alpha}$;
        \item[2.] $U_{\varphi(\alpha)}\subset
        DA(\varphi(\alpha'))$, for any $\alpha'\in U_{\alpha}$
    \end{itemize}
\end{corollary}

\begin{proof}
For $\alpha\in D_{1}$ and $x\in DA(\varphi(\alpha))$, the function
$V(x,\alpha)$ from Theorem 2 is considered. The real and
non-negative function $V$ is defined on the open set
$G=\bigcup\limits_{\alpha\in D_{1}} (DA(\varphi(\alpha))\times
\{\alpha\}) \subset \Omega\times D_{1}$, it is continuous and
equal to zero on the set $\Gamma
=\{(\varphi(\alpha),\alpha)/\alpha\in D_{1}\}\subset G$.

As $V$ is continuous and it is equal to zero in
$(\varphi(\alpha),\alpha)\in G$, there is an open neighborhood $W$
of $(\varphi(\alpha),\alpha)$ such that for any $(x',\alpha')\in
W$, the inequality $V(x',\alpha')<1$ holds. Let be $U_{\alpha}$ an
open neighborhood of $\alpha$ and $U_{\varphi(\alpha)}$ of
$\varphi(\alpha)$ such that $U_{\varphi(\alpha)}\times
U_{\alpha}\subset W$. As the function $\varphi$ is continuous, it
can be admitted that for any $\alpha'\in U_{\alpha}$, we have
$\varphi(\alpha')\in U_{\varphi(\alpha)}$ (contrarily, the
neighborhood $U_{\alpha}$ can be replaced with a smaller
neighborhood $U_{\alpha}'\subset U_{\alpha}$, for which we have
$\varphi(\alpha')\in U_{\varphi(\alpha)}$, for any $\alpha'\in
U_{\alpha}'$).

Thus, for any $(x',\alpha')\in U_{\varphi(\alpha)}\times
U_{\alpha}$, we have $V(x',\alpha')<1$. This means that for any
$x'\in U_{\varphi(\alpha)}$ and any $\alpha'\in U_{\alpha}$, we
have that $x'\in DA(\varphi(\alpha'))$. Thus,
$U_{\varphi(\alpha)}\subset DA(\varphi(\alpha'))$, for any
$\alpha'\in U_{\alpha}$.
\end{proof}

\begin{rem}
 Corollary 1 states that for any $\alpha'\in U_{\alpha}$, both
 maneuvers $\alpha\rightarrow\alpha'$ and
 $\alpha'\rightarrow\alpha$ are successful on the path $\varphi$.
\end{rem}

\begin{thm}
For two steady states $\varphi(\alpha^{\star})$ and
$\varphi(\alpha^{\star\star})$ belonging to the analytic path
$\varphi$ of asymptotically stable steady states of (\ref{ec}),
there exist a finite number of values of the control parameters
$\alpha^{1},\alpha^{2},...,\alpha^{p}\in D_{1}$ such that all the
maneuvers
\begin{equation}\label{mansuccesive}
    \alpha^{\star}\rightarrow\alpha^{1}\rightarrow\alpha^{2}\rightarrow
    ...\rightarrow\alpha^{p}\rightarrow\alpha^{\star\star}
\end{equation}
are successful on the path $\varphi$.
\end{thm}

\begin{proof}
Let be $P\subset D_{1}$ a polygonal line which joins
$\alpha^{\star}$ and $\alpha^{\star\star}$. For any $\alpha\in P$
we consider the neighborhoods $U_{\alpha}$ and
$U_{\varphi(\alpha)}$ given by Corollary 1.

The family of neighborhoods $\{U_{\alpha}\}_{\alpha\in P}$ is a
covering with open sets of the compact polygonal line $P$. From
this covering we can subtract a finite covering of $P$, i.e.,
there exist
$\bar{\alpha}^{1},\bar{\alpha}^{2},...,\bar{\alpha}^{q}\in P$ such
that $P\subset\bigcup\limits_{k=1}^{q}U_{\bar{\alpha}^{k}}$. More,
it can be assumed that $\alpha^{\star}\in U_{\bar{\alpha}^{1}}$
and $\alpha^{\star\star}\in U_{\bar{\alpha}^{q}}$ and that the
intersections $U_{\bar{\alpha}^{k}}\cap P$ are open and connected
sets in $P$, and
\begin{equation}
\nonumber
    (U_{\bar{\alpha}^{k}}\cap P)\cap(U_{\bar{\alpha}^{k+2}}\cap
    P)=\emptyset \qquad \textrm{for any } k=1,2,...,q-2.
\end{equation}

Taking into account Remark 2, as $\alpha^{\star}\in
U_{\bar{\alpha}^{1}}$ and $\alpha^{\star\star}\in
U_{\bar{\alpha}^{q}}$, it comes naturally that the maneuvers
$\alpha^{\star}\rightarrow\bar{\alpha}^{1}$ and
$\bar{\alpha}^{q}\rightarrow\alpha^{\star\star}$ are successful on
the path $\varphi$.

We still have to prove that each maneuver
$\bar{\alpha}^{k}\rightarrow\bar{\alpha}^{k+1}$ is successful for
any $k=1,2,...,q-1$.

If $\bar{\alpha}^{k}\in U_{\bar{\alpha}^{k+1}}$, Remark 2 provides
that the maneuver $\bar{\alpha}^{k}\rightarrow\bar{\alpha}^{k+1}$
is successful on the path $\varphi$.

If $\bar{\alpha}^{k}\notin U_{\bar{\alpha}^{k+1}}$, a point
$\bar{\alpha}^{k,k+1}\in (U_{\bar{\alpha}^{k}}\cap P)\cap
(U_{\bar{\alpha}^{k+1}}\cap P)$ is considered. Remark 2 provides
that both maneuvers
$\bar{\alpha}^{k}\rightarrow\bar{\alpha}^{k,k+1}$ and
$\bar{\alpha}^{k,k+1}\rightarrow\bar{\alpha}^{k+1}$ are successful
on the path $\varphi$.

Thus, eventually inserting control parameters
$\bar{\alpha}^{k,k+1}$ between $\bar{\alpha}^{k}$ and
$\bar{\alpha}^{k+1}$, we come to find (after changing the notation
and re-numbering) a finite sequence
$\alpha^{1},\alpha^{2},...,\alpha^{p}\in D_{1}$ such that all the
maneuvers
\begin{equation}\nonumber
    \alpha^{\star}\rightarrow\alpha^{1}\rightarrow\alpha^{2}\rightarrow
    ...\rightarrow\alpha^{p}\rightarrow\alpha^{\star\star}
\end{equation}
are successful on the path $\varphi$.
\end{proof}

\begin{rem}
Theorem 3 states that two steady states belonging to an analytic
path of asymptotically stable steady states can be transferred one
in the other using a finite number of successful maneuvers along
the considered path.
\end{rem}

\section{Numerical examples}

\subsection*{Example 1.}

The following one-dimensional discrete dynamical system with
control is considered:
\begin{equation}\label{ex1}
    x_{k+1}=\alpha x_{k}(1-x_{k}), \qquad k=0,1,2,...
\end{equation}
where $x\in\mathbb{R}$ is the state parameter and
$\alpha\in\mathbb{R}$ is the control parameter. This dynamical
system is frequently used for showing chaotic behavior and it is
subject to several themes of research (\cite{Feigenbaum}). We will
illustrate using this system the concepts of analytic path,
analytic path of asymptotically stable steady states, domain of
attraction, successful maneuver.

For $\alpha\neq 0$, the steady states of (\ref{ex1}) are $x=0$ and
$x=\frac{\alpha-1}{\alpha}$, while for $\alpha=0$ corresponds only
the $x=0$ steady state. Thus, for (\ref{ex1}) there are three
paths of steady states:
$\varphi_{1}(\alpha)=\frac{\alpha-1}{\alpha}$, for $\alpha<0$,
$\varphi_{2}(\alpha)=0$, for $\alpha\in\mathbb{R}$ and
$\varphi_{3}(\alpha)=\frac{\alpha-1}{\alpha}$, for $\alpha>0$,.
The path $\varphi_{1}$ contains only unstable steady states. The
steady states belonging to the path $\varphi_{2}$ are
asymptotically stable for $\alpha\in (-1,1)$. The steady states
belonging to the path $\varphi_{3}$ are asymptotically stable for
$\alpha\in (1,3)$.

In Fig. \ref{fig-ex1-ram}, the three paths of steady states are
plotted. In Fig. \ref{fig-ex1-dom} the gray rectangle represents
the reunion of the domains of attraction of the asymptotically
stable steady states of $\varphi_{3}$, while the vertical gray
line denotes the domain of attraction of the steady state
$\varphi_{2}(0)$. In both figures, the black parts of the paths of
steady states represent the asymptotically stable steady states
while the gray parts of the paths represent the unstable steady
states.

\renewcommand{\thefigure}{\arabic{figure}.\arabic{subfigure}}
\setcounter{subfigure}{1}

\begin{figure}[htbp]
\begin{minipage}[t]{0.5\linewidth}
\centering
\includegraphics*[bb=3cm 0cm 13.5cm
6.5cm, width=6.3cm]{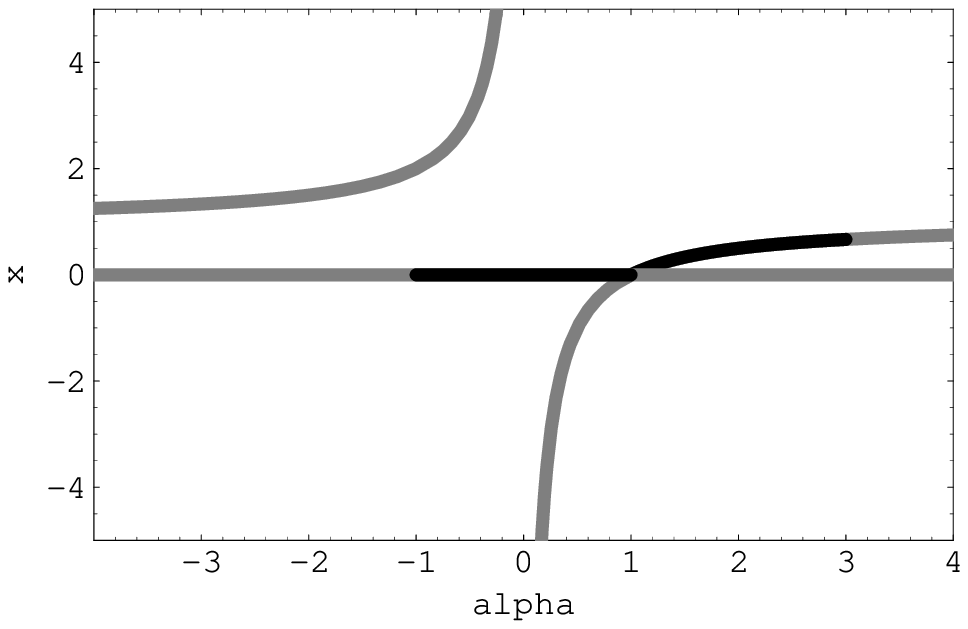} \caption{The paths of steady states
for (\ref{ex1})}\label{fig-ex1-ram}
\end{minipage}
\addtocounter{figure}{-1} \addtocounter{subfigure}{1}
\begin{minipage}[t]{0.5\linewidth}
\centering
\includegraphics*[bb=3cm 0cm 13.5cm
6.5cm, width=6.5cm]{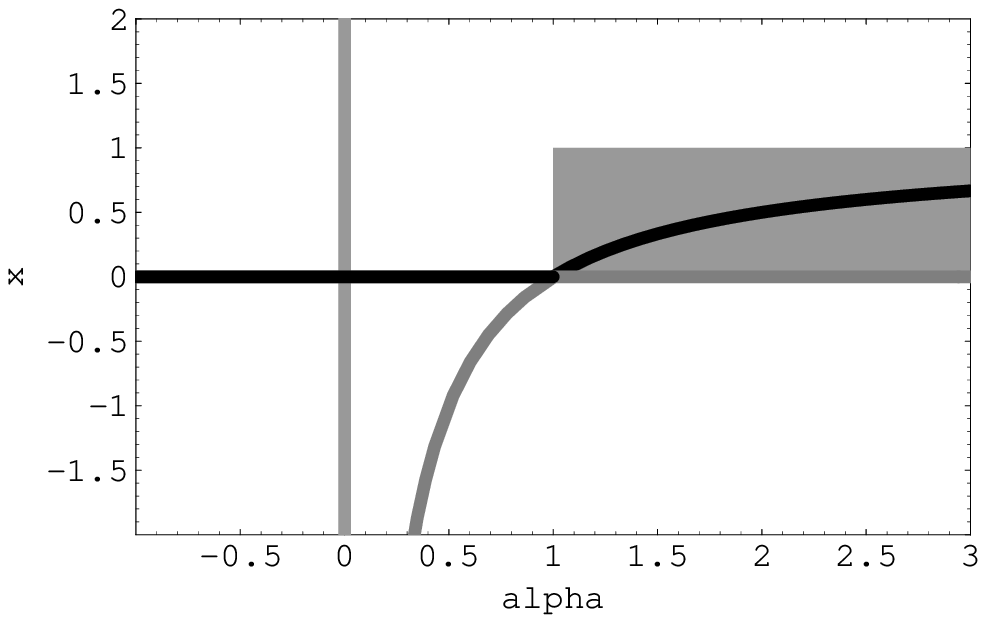} \caption{The domains of attractions
of the asymptotically stable steady states of the path
$\varphi_{3}$ and for the asymptotically stable steady state
$\varphi_{2}(0)=0$ of $\varphi_{2}$}\label{fig-ex1-dom}
\end{minipage}
\end{figure}

The domain of attraction of the steady state $\varphi_{2}(0)=0$ is
$DA(\varphi_{2}(0))=\mathbb{R}$, while the domain of attraction of
a steady state $\varphi_{3}(\alpha)$ for $\alpha\in(1,3)$ is
$DA(\varphi_{3}(\alpha))=(0,1)$. These domains of attraction can
be obtained using the staircase method \cite{Kel-Pet}, or can be
estimated numerically using the method described in \cite{EBBB}).

The steady state $\varphi_{3}(1.1)$ can be directly transferred by
a single maneuver to $\varphi_{3}(2.9)$, because
$\varphi_{3}(1.1)\in DA(\varphi_{3}(2.9))=(0,1)$. The
$DA(\varphi_{3}(1.1))$ includes $\varphi_{3}(2.9)$, thus, the
maneuver $\alpha:2.9\rightarrow 1.1$ is also successful.

All asymptotically stable steady states of $\varphi_{3}$ are in
the domain of attraction of the steady state $\varphi_{2}(0)$.
This means that every maneuver $\alpha\rightarrow 0$, for
$\alpha\in(1,3)$ is successful between the paths $\varphi_{3}$ and
$\varphi_{2}$. Though, a steady state $\varphi_{2}(\alpha)$ cannot
be transferred in an asymptotically stable steady state of
$\varphi_{3}$, because any maneuver of the type
$\alpha\rightarrow\alpha'$, with $\alpha'\in (1,3)$ causes a
transfer to the unstable steady state $\varphi_{2}(\alpha')=0$.

\subsection*{Example 2.}

The following one-dimensional discrete dynamical system with
control is considered:
\begin{equation}\label{ex2}
    x_{k+1}=(x_{k}-\alpha)^{3}+\alpha, \qquad k=0,1,2,...
\end{equation}
where $x\in\mathbb{R}$ is the state parameter and
$\alpha\in\mathbb{R}$ is the control parameter.

The sequence $x_{k}$, with the starting point $x_{0}$ which
satisfies (\ref{ex2}) is:

\begin{equation}\label{ex2-formula}
    x_{k}=(x_{0}-\alpha)^{3^{k}}+\alpha, \qquad k=0,1,2,...
\end{equation}

There are three analytic paths of steady states for (\ref{ex2}):
$\varphi_{1}(\alpha)=\alpha$, $\varphi_{2}(\alpha)=\alpha-1$ and
$\varphi_{3}(\alpha)=\alpha+1$, defined for $\alpha\in\mathbb{R}$.
The path $\varphi_{1}$ is an analytic path of asymptotically
stable steady states while $\varphi_{2}$ and $\varphi_{3}$ are
analytic paths of unstable steady states. In Fig.
\ref{fig-ex2-ram}, the continuous line represents the path
$\varphi_{1}$, while the dashed lines represent the paths
$\varphi_{2}$ and $\varphi_{3}$.

For any $\alpha\in\mathbb{R}$, the domain of attraction of the
asymptotically stable steady state $\varphi_{1}(\alpha)$ is
$DA(\varphi_{1}(\alpha))=(\alpha-1,\alpha+1)$.

For $\alpha^{\star}=0$ and $\alpha^{\star\star}=2$, let's consider
the asymptotically stable steady states
$\varphi_{1}(\alpha^{\star})=0$ and
$\varphi_{1}(\alpha^{\star\star})=2$. The maneuver
$\alpha:\alpha^{\star}=0\rightarrow 2=\alpha^{\star\star}$ is not
successful, because $\varphi_{1}(\alpha^{\star})=0\notin
DA(\varphi_{1}(\alpha^{\star\star})=2)=(1,3)$. Though, a finite
number of maneuvers can be found, which transfer the steady state
$\varphi_{1}(\alpha^{\star})=0$ to the steady state
$\varphi_{1}(\alpha^{\star\star})=2$, for example:

\begin{equation}\label{ex2-man}
    \alpha: \alpha^{\star}=0\rightarrow 0.7 \rightarrow 1.4
    \rightarrow 2=\alpha^{\star\star}
\end{equation}

\renewcommand{\thefigure}{\arabic{figure}}
\setcounter{figure}{1}

\begin{figure}[htbp]
\centering
\includegraphics*[bb=3cm 0cm 13.5cm
6.5cm, width=7cm]{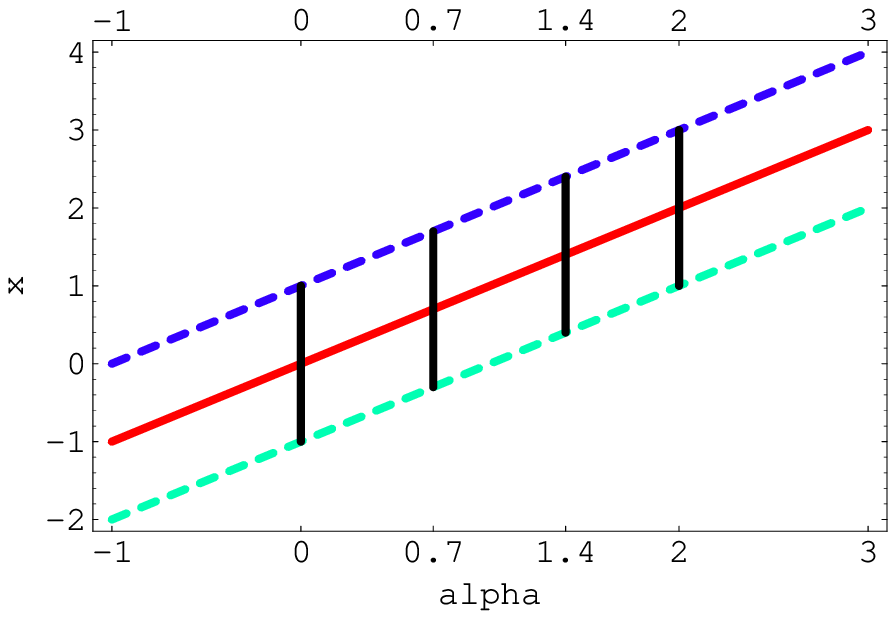} \caption{The paths of steady states
for (\ref{ex2}) and the domains of attraction of the steady states
of $\varphi_{1}$ corresponding to the control parameters used in
\ref{ex2-man}}\label{fig-ex2-ram}
\end{figure}

These maneuvers are successful, because
$\varphi_{1}(\alpha^{\star})\in DA(\varphi_{1}(0.7))=(-0.3,1.7)$,
$\varphi_{1}(0.7)\in DA(\varphi_{1}(1.4))=(0.4,2.4)$ and
$\varphi_{1}(1.4)\in DA(\varphi_{1}(\alpha^{\star\star}))=(1,2)$.
In Fig. \ref{fig-ex2-ram}, the vertical segments represent the
domains of attraction of the steady states corresponding to the
maneuvers (\ref{ex2-man}).

\subsection*{Example 3.}

The following two-dimensional discrete dynamical system with
control is considered:
\begin{equation}\label{ex3}
\left\{%
\begin{array}{ll}
    x_{k+1}=(x_{k}-\alpha)[(x_{k}-\alpha)^{2}+(y_{k}-\alpha)^2]+\alpha \\
    y_{k+1}=(y_{k}-\alpha)[(x_{k}-\alpha)^{2}+(y_{k}-\alpha)^2]+\alpha \\
\end{array}%
\right.
\end{equation}
where $(x,y)\in\mathbb{R}^{2}$ is the state parameter and
$\alpha\in\mathbb{R}$ is the control parameter.

There are an infinity of analytic paths of steady states for
(\ref{ex3}): $\varphi(\alpha)=(\alpha,\alpha)$ and
$\varphi_{t}(\alpha)=(\alpha+\cos t,\alpha+\sin t)$, for $t\in
[0,2\pi)$; all paths are defined for $\alpha\in\mathbb{R}$. The
path $\varphi$ is an analytic path of asymptotically stable steady
states while $\varphi_{t}$ are analytic paths of unstable steady
states, for any $t\in [0,2\pi)$.

For any $\alpha\in\mathbb{R}$, the domain of attraction of the
asymptotically stable steady state
$\varphi(\alpha)=(\alpha,\alpha)$ is the ball
$B((\alpha,\alpha),1)=\{(x,y)\in\mathbb{R}^2/(x-\alpha)^2+(y-\alpha)^2<1\}$.

\begin{figure}[htbp]
\centering
\includegraphics*[bb=3cm 0cm 13cm
10.5cm, width=6cm]{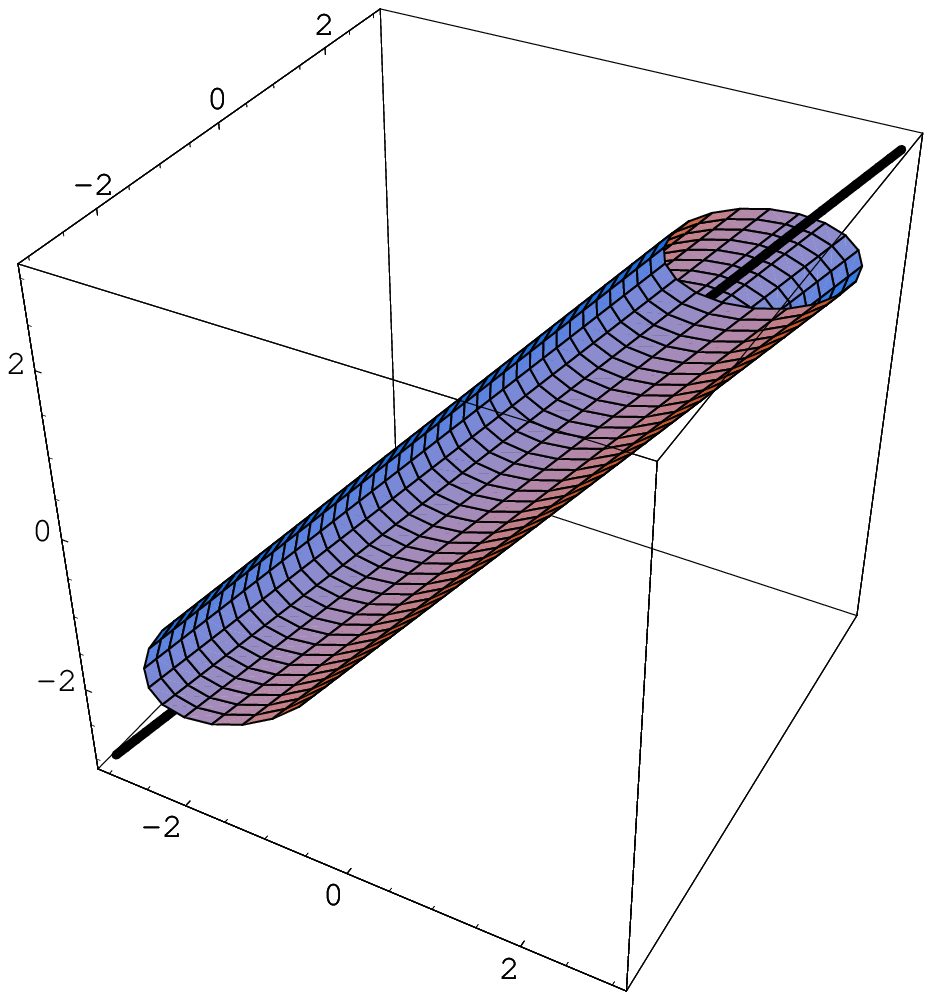} \caption{The paths of steady states
for (\ref{ex3})}\label{fig-ex2-ram}
\end{figure}

For $\alpha^{\star}=-1$ and $\alpha^{\star\star}=1$, let's
consider the asymptotically stable steady states
$\varphi(\alpha^{\star})=(-1,-1)$ and
$\varphi(\alpha^{\star\star})=(1,1)$. The maneuver
$\alpha:\alpha^{\star}=-1\rightarrow 1=\alpha^{\star\star}$ is not
successful, because $\varphi(\alpha^{\star})=(-1,-1)\notin
DA(\varphi(\alpha^{\star\star})=(1,1))=B((1,1),1)$. Though, a
finite number of maneuvers can be found, which transfer the steady
state $\varphi(\alpha^{\star})=(-1,-1)$ to the steady state
$\varphi(\alpha^{\star\star})=(1,1)$, for example:

\begin{equation}\label{ex3-man}
    \alpha: \alpha^{\star}=-1\rightarrow -0.5 \rightarrow 0
    \rightarrow 0.5\rightarrow 1=\alpha^{\star\star}
\end{equation}

These maneuvers are successful, because
$\varphi(\alpha^{\star})\in DA(\varphi(-0.5))=B((-0.5,-0.5),1)$,
$\varphi(-0.5)\in DA(\varphi(0))=B((0,0),1)$, $\varphi(0)\in
DA(\varphi(0.5))=B((0.5,0.5),1)$, and $\varphi(0.5)\in
DA(\varphi(\alpha^{\star\star}))=B((1,1),1)$.

%%%%%%%%%%%%%%%%%%%%%%%%%%%%%%%%%%%%%%%%%%%%%%%%%%%%%%%%%%%%%%%%%%%%

\end{document}